\documentclass{amsart}

\usepackage[english]{my-shortcuts}
\usepackage{scalerel,stackengine}

\usepackage{a4wide}
\usepackage{amsmath,amsfonts,amssymb,latexsym,amsthm}
\usepackage{dsfont}
\usepackage{comment} 
\usepackage{graphicx}
\usepackage{caption} 

\newcommand{\ZZ}{\mathbb{Z}}

\stackMath
\newcommand\reallywidehat[1]{%
\savestack{\tmpbox}{\stretchto{%
  \scaleto{%
    \scalerel*[\widthof{\ensuremath{#1}}]{\kern-.6pt\bigwedge\kern-.6pt}%
    {\rule[-\textheight/2]{1ex}{\textheight}}
  }{\textheight}%
}{0.5ex}}%
\stackon[1pt]{#1}{\tmpbox}%
}
\begin{document}

\title[Closed forms related to integrals of $\zeta$ and $L$-functions]{
Weighted mean squares of $\zeta$ and $L(\cdot,\chi_4)$: closed forms,
reciprocity, and a missing identity}
\author{S\'ebastien Darses -- Erwan Hillion}
\date{}


\address{PyxiScience\\ \ CNRS -- Universit\'e de Montr\'eal CRM–CNRS}
\email{seb.darses@gmail.com}

\address{Acad\'emie Militaire Saint-Cyr Co\"etquidan -- CREC}
\email{erwan.hillion@st-cyr.terre-net.defense.gouv.fr}

\maketitle

\begin{abstract}


This short note is a concise yet complete overview, intended for students, researchers, and AIs, bringing together in a unified notation the closed-form and reciprocity identities arising from weighted mean squares of the Riemann zeta function and Dirichlet $L$-functions on the critical line.

In addition, we prove a new reciprocity identity: the exponentially weighted counterpart of a formula of Lewis and Zagier for the primitive character of modulus $4$.
\end{abstract}

\section{Introduction}

Let $\zeta$ be the Riemann zeta function and $L(\cdot,\chi)$ be the Dirichlet $L$-function associated to a Dirichlet character $\chi$.
We consider here the following weights on the critical line $s=\frac{1}{2}+it$, $t\in\R$:
$$\frac{1}{s(1-s)}=\frac{1}{t^2+\frac{1}{4}}, \quad \frac{\Gamma(s)\Gamma(1-s)}{\pi}=\frac{1}{\cosh(\pi t)},$$
which we call the $sq-$weight and $ch-$weight respectively, where $\Gamma$ is the usual Gamma function.

Let $\chi_4$ be the primitive odd character of conductor $4$. It is the real-valued $4-$periodic function defined by $\chi_4(0)=0$, $\chi_4(1)=1$, $\chi_4(2)=0$, and $\chi_4(3)=\chi_4(-1)=-1$.
Let $\{x\}=x-\lfloor x \rfloor$ denote the fractional part of a real number $x\geq 0$, where $\lfloor x \rfloor$ is its integer part. 


We first consider the following absolutely convergent integrals (e.g., the first one is called a twisted mean-square in \cite{BC13a}) for $\alpha>0$:
\begin{eqnarray*}
I_{\zeta,sq}(\alpha) = \int_{-\infty}^\infty \left|\zeta\left(\tfrac{1}{2}+it\right)\right|^2 \alpha^{it} \frac{dt}{t^2+\tfrac{1}{4}}, &  &
I_{\zeta,ch}(\alpha) = \int_{-\infty}^\infty \left|\zeta\left(\tfrac{1}{2}+it\right)\right|^2 \alpha^{it} \frac{dt}{\cosh(\pi t)}, \\
    I_{\chi_4,sq}(\alpha) = \int_{-\infty}^\infty \left|L\left(\tfrac{1}{2}+it,\chi_4\right)\right|^2 \alpha^{it} \frac{dt}{t^2+\tfrac{1}{4}}, & &
    I_{\chi_4,ch}(\alpha) = \int_{-\infty}^\infty \left|L\left(\tfrac{1}{2}+it,\chi_4\right)\right|^2 \alpha^{it} \frac{dt}{\cosh(\pi t)}.
\end{eqnarray*}

They all satisfy the following remarkable {\em closed-form} property:

\begin{theo}
Let $\alpha=n/m>0$ be a rational number where $n,m$ are coprime. Then $I_{\zeta,sq}(\alpha)$, $I_{\zeta,ch}(\alpha)$, $I_{\chi_4,sq}(\alpha)$, $I_{\chi_4,ch}(\alpha)$ can be written as {\em finite sums,
whose lengths are given by $n$ and $m$}, involving the functions $\{\cdot\},\log, \cot$, together with Euler's constant $\gamma$, $\pi$, and rational coefficients.
\end{theo}

The underlying sums are cotangent sums, and display a symmetry under $n\leftrightarrow m $, a reciprocity phenomenon discussed in Section 4.
The first three closed-form formulas are known:

\begin{enumerate}
    \item 
The first identity, concerning $I_{\zeta,sq}(n/m)$, is Vasyunin's formula \cite{Vas95}. The reformulation with Mellin isometry, considered here, is explicitly written in \cite{BC13a} p. 225.
\item 
The second one, concerning $I_{\zeta,ch}(n/m)$, is a consequence of a reciprocity formula of Bettin and Conrey in \cite{BC13b}, and is also proved by other means in \cite{DH21a}. 
\item 
The third one, concerning $I_{\chi_4,sq}(n/m)$, is a formula of Lewis and Zagier in \cite{LZ19}. 
\item The last one, concerning $I_{\chi_4,ch}(n/m)$, seems to be new and is proved here, which happily completes the above rectangle. 
\end{enumerate}

In Section \ref{tools}, we develop the tools needed to derive these identities, which we then collect and discuss in Section \ref{review}. Section \ref{proof} is devoted to the proof of the new closed-form formula for $I_{\chi_4,ch}(n/m)$.
Finally, in Section \ref{review-moment}, we put in a unified notation the closed forms related to
\begin{eqnarray*}
\int_{-\infty}^\infty  t^{N}  \left|\zeta\left(\tfrac{1}{2}+it\right)\right|^2 \frac{dt}{\cosh(\pi t)}, & \quad & 
\int_{-\infty}^\infty t^{N} \left|L\left(\tfrac{1}{2}+it, \chi\right)\right|^{2} \frac{dt}{\cosh(\pi t)},
\end{eqnarray*}
the first one obtained in \cite{DH24} and the second one in \cite{DRR26} (D., Ringeling and Royer) for any primitive Dirichlet character $\chi$.

Let us notice the important, though elementary, fact that all these closed forms 
provide exact characterizations of $|\zeta|$ and $\left|L\left(\cdot, \chi_4\right)\right|$ on the critical line by means of Fourier transforms. Indeed, on the one hand, $\log(\Q_{>0})$ is dense in $\R$ and the Fourier transforms, associated to the twisted mean-squares, are both continuous and injective (on $L^1(\R)$ e.g.). On the other hand, the corresponding moment problems are determinate, as seen in \cite{DH24} and \cite{DRR26}; these moments are the $N$th derivatives at $0$ of the corresponding Fourier transforms. 

As the reader will notice in the sequel, none of these identities involve the divisor function $d(n)$ (the number of divisors of an integer $n$), while this arithmetic function is as the core definition of the Dirichlet series (see e.g. \cite{Tit86,Ten22})
$$
\zeta(s)^2 = \sum_{a,b\geq1}\frac{1}{a^s b^s}=\sum_{n= 1}^\infty\frac{d(n)}{n^s}, \quad L(s,\chi)^2 = \sum_{n= 1}^\infty\frac{\chi(n)d(n)}{n^s}, \qquad \mathfrak{R}(s)>1.
$$

The interest in closed forms dates back to the Basel problem, posed by Mengoli and solved by Euler: $\zeta(2)=\pi^2/6$, and have since concerned many more aspects of values of $\zeta$, which is a huge topic in number theory, connected to Bernoulli numbers, Stirling numbers, partition functions, polylogarithm, polyzetas among many others objects. See, e.g., remarkable and delightful expositions in \cite{ORS17}, \cite{AIKZ14}, and \cite{GKP94}.

Regarding many insightful discussions on closed forms we refer to the beautiful paper \cite{BoC13}. 


\section{Special simple forms}

As an {\em amuse-bouche} for the curious reader, let us first mention a few remarkable identities without parameters:
\begin{eqnarray*}
\int_{-\infty}^\infty \zeta\left(\tfrac{1}{2}+it\right) \frac{dt}{t^2+\tfrac{1}{4}}  & = & 2\pi (\gamma - 1) \ =\ -2.6564323\ldots\\
\int_{-\infty}^\infty \left|\zeta\left(\tfrac{1}{2}+it\right)\right|^2 \frac{dt}{t^2+\tfrac{1}{4}}  & = & 2\pi (\log(2\pi)-\gamma) \ =\  7.920969\ldots \\
\int_{-\infty}^\infty \left|\zeta\left(\tfrac{1}{2}+it\right)\right|^2 \frac{dt}{\cosh(\pi t)}  & = & 2 (\log(2\pi)-\gamma) - 1 \ =\  1.521322\ldots \\
\int_{-\infty}^\infty \left|\zeta\left(\tfrac{1}{2}+it\right)\right|^2\  \frac{3-\sqrt{8}\cos(t\log 2)}{t^2+\tfrac{1}{4}}dt  & = & 2\pi \log 2 \ =\  4.3551721\ldots 
\end{eqnarray*}

We like how $-1,\pi,\gamma$ and $\log$ dance, appear and disappear! (Look carefully.) 

We noticed the first identity on the website Mathoverflow \cite{Math18}.  
The sign is consistent with the (small) concentration within the integral around the value $\zeta(\tfrac {1}{2})=-1.4603545\ldots$ 
No closed form of $\zeta(\tfrac {1}{2})$ is currently known, while one knows: $\displaystyle \Gamma \left({\tfrac {1}{2}}\right)={\sqrt {\pi }}$, and $\d\frac{\zeta'(\tfrac {1}{2})}{\zeta(\tfrac {1}{2})}= \frac{\log(8\pi)}{2} +\frac{\pi}{4} + \frac{\g}{2}$.

The second identity appears in \cite{Cof11}, and can be obtained taking $n=m=1$ in Theorem \ref{th:Vasyunin}.

The third one can be obtained taking $n=m=1$ in Theorem \ref{DH21}, or $N=0$ in Theorem \ref{DH24}.

The last one is due to Ivi\'c \cite{Ivi03}
and stems from his beautiful identity:
\begin{eqnarray}
\int_{-\infty}^{\infty}
\left| \frac{1-2^{1-s}}{s} \zeta(s)\right|^2 \frac{dt}{2\pi}
& =  &
\frac{1-2^{1-2\sg}}{2\sg}\ \zeta(2\sg),
\qquad s=\sigma+it,
\end{eqnarray} 
which he first computed for $\sg>1$, and then obtained by analytic continuation for $\sg>0$. One only needs to take $\sg\to 1/2$ in both sides. The disappearance of $\gamma$ in the last expression is explained by the absence of the pole of $\zeta$ in $(1-2^{1-s})\zeta(s)$.

Integrals involving $\zeta$ or the $\xi$-function, their relations with special functions and associated symmetries, are an endless story: we might start with \cite{Ram15}. See, e.g., \cite{BD10, DK21, DSS24, Mil24} and numerous references therein for recent insights.

Regarding $L(\cdot,\chi_4)$, we have the remarkable values:
\begin{eqnarray*}
\int_{-\infty}^\infty \left|L\left(\tfrac{1}{2}+it,\chi_4\right)\right|^2 \frac{dt}{t^2+\tfrac{1}{4}}  & = & \frac{\pi^2}{2} \ =\  4.934802\ldots \\
\int_{-\infty}^\infty \left|L\left(\tfrac{1}{2}+it,\chi_4\right)\right|^2 \frac{dt}{\cosh(\pi t)}  & = & \frac{1}{2}.
\end{eqnarray*}
The first identity can be obtained, for instance, by taking $n=m=1$ in Theorem \ref{LZ19}. The second can be deduced from Theorem \ref{DRR26} by taking $N=0$ and $\chi=\chi_4$, but it is much easier to get it directly with Lemma \ref{mellin-iso} and
$\d
\int_0^\infty \frac{dx}{\cosh^2 x}
= \left[\tanh x\right]_0^\infty
= 1.
$

\bigskip

\section{Tools to analyse the twisted mean-squares of $\zeta$ and $L(\cdot,\chi_4)$} \label{tools}

\subsection{Classical Mellin transforms}
Let 
\begin{eqnarray*}
    f(x) = \sum_{n\geq 1} e^{-nx} - \frac{1}{x} = \frac{1}{e^{x}-1}-\frac{1}{x}, \quad x>0.
\end{eqnarray*}
The following Mellin transform, for $0<\mathfrak{R}(s)<1$,
\begin{eqnarray} \label{mellin-gamma-zeta}
\cal M[f](s) := \int_0^\infty f(x) x^{s-1}dx = \int_0^\infty \left(\frac{1}{e^{x}-1}-\frac{1}{x}\right) x^{s-1}dx = \Gamma(s)\zeta(s),
\end{eqnarray}
can be generalized to $L-$functions. We also write $\sg= \mathfrak{R}(s)$.

Let $\chi$ be a non principal Dirichlet character of modulus $q$. We have 
the following generalization of $f$ for all $x>0$:
\begin{eqnarray*}
    f_\chi(x) = \sum_{n\geq 1} \chi(n) e^{-nx} 
            =  \sum_{m\geq 0}\sum_{k=1}^{q-1} \chi(mq+k) e^{-(mq+k)x} 
         =  \sum_{k=1}^{q-1} \chi(k) e^{-kx}\sum_{m\geq 0} e^{-mq x} 
             =  \frac{F_{\chi}(e^{-x})}{1-e^{-qx}},
\end{eqnarray*}
where $\d F_\chi(X)=\sum_{k=1}^{q-1} \chi(k) X^k$ is a Fekete polynomial. See \cite{CGPS00} for many important properties of this object. 

The function $f_\chi$ is regular at $0$ since $X-1\ |\ F_\chi(X)$, and  $f_\chi$ is exponentially decreasing at infinity since $X\ |\ F_\chi(X)$. We then have the classical representation (see e.g. \cite{IR90})
\begin{eqnarray} \label{mellin-g-L}
\cal M[f_\chi](s) & = & \Gamma(s) L(s,\chi), \quad \mathfrak{R}(s)> 0.
\end{eqnarray}
Recall that $L(\cdot,\chi)$ has no pole at $1$, and thus there is no additional factor $1/x$ as in (\ref{mellin-gamma-zeta}).

We consider here one of the "simplest" $L$-functions:
\begin{eqnarray}
L(s,\chi_4) & = & 1-\frac{1}{3^s}+\frac{1}{5^s}- \frac{1}{7^s} +\cdots
\end{eqnarray}
for $\mathfrak{R}(s)> 1$ and then continued analytically.
We have: 
\begin{eqnarray}
\int_0^\infty \frac{x^{s-1}}{\cosh(x)} dx & = & 2 \ \Gamma(s) L(s,\chi_4),
\end{eqnarray}
since
\begin{eqnarray*}
    f_{\chi_4}(x) = \frac{F_{{\chi_4}}(e^{-x})}{1-e^{-qx}} = \frac{e^{-x}-e^{-3x}}{1-e^{-4x}} = \frac{e^{x}-e^{-x}}{e^{2x}-e^{-2x}} 
        = \frac{1}{e^{x} + e^{-x}}.
\end{eqnarray*}
A remarkable feature of $f_{\chi_4}$ is that it is an even special function.

Let $\mathds{1}_{A}$ be the indicator function of a set $A$. 

Following Lewis and Zagier \cite{LZ19}, let us set (see \cite[(12) p.5]{LZ19}):
\begin{eqnarray}
    S(x) = {\sum_{0<k\leq x}}^* \chi(k) := \sum_k \mathds{1}_{(4k+1,4k+3)}(x) + \frac{1}{2}\sum_{k} \mathds{1}_{2k+1}(x),
\end{eqnarray}
which is a generalization of $\lfloor x \rfloor = \sum_{k\leq x}1$.
They note that (see \cite[(21) p.9]{LZ19}):
\begin{eqnarray}
    \int_0^\infty S\left(\frac{1}{x}\right) x^{s-1}dx & = & \frac{L(s,\chi_4)}{s}, \quad \sg>0,
\end{eqnarray}
to be compared with:
\begin{eqnarray}
    \int_0^\infty \left\lfloor\frac{1}{x}\right\rfloor x^{s-1}dx & = & \frac{\zeta(s)}{s}, \quad \sg>1, \nonumber\\
    \int_0^\infty \left\{\frac{1}{x}\right\} x^{s-1}dx & = & -\frac{\zeta(s)}{s}, \quad 0<\sg<1.
\end{eqnarray}

The change of variable $x\to nx$ gives the following expressions : $\displaystyle \int_0^\infty \left\{\frac{1}{nx}\right\} x^{s-1}dx = -n^{-s}\frac{\zeta(s)}{s}$, $\d \int_0^\infty \left(\frac{1}{mt}-\frac{1}{e^{mt}-1}\right)t^{s-1}dt = m^{-s} \Gamma\zeta(s)$, etc., which will be used in the Mellin isometries below.


\subsection{Auto-correlation functions and Mellin isometry}

As done in \cite[p.225]{BC13a} for Vasyunin's formula, where the authors use the Mellin isometry:
\begin{eqnarray*}
    \int_{0}^\infty \left\{\frac{1}{nx}\right\}\left\{\frac{1}{mx}\right\} dx & = &
    \frac{1}{2\pi i}\int_{-\infty}^\infty n^{-s} m^{-\b s} \left|\frac{\zeta(s)}{s}\right|^2 ds, \qquad s=\frac{1}{2}+it,
\end{eqnarray*}
we can then relate all the integrals under study to the corresponding multiplicative autocorrelation functions for $\lb>0$ and $\mathfrak{R}(z)>0$:
\begin{eqnarray*}
    A_{\zeta, sq}(\lb) = \int_0^\infty \left\{\frac{1}{x}\right\}\left\{\frac{1}{\lb x}\right\} dx, & \quad &
    A_{\zeta, ch}(z) = \int_0^\infty \left(\frac{1}{x}-\frac{1}{e^{x}-1}\right)\left(\frac{1}{zx}-\frac{1}{e^{zx}-1}\right) dx, \\
    A_{\chi_4,sq}(\lb) = \int_0^\infty S\left(\frac{1}{x}\right) S\left(\frac{1}{\lb x}\right) dx, & \quad  &
    A_{\chi_4,ch}(z) = \int_{0}^\infty \frac{dx}{\cosh(x)\cosh(z x)}.
\end{eqnarray*}

The word "autocorrelation function" has been introduced in this context in \cite{BDBLS05, BM18} where the authors provide a detailed study of the function: $\R_{\geq0}\to \R_{>0},\ \lb\mapsto A_{\zeta, sq}(1/\lb)$. 

Considering the $ch-$weight instead of the $sq-$weight induces a huge regularization effect on the corresponding autocorrelation functions: $A_{\zeta, ch}$ and $A_{\chi_4,ch}$ have an analytic continuation to $\C\setminus \R_{\leq0}$ (cf. \cite{DN24, DRR26}) while $A_{\zeta, sq}$ is continuous but non-differentiable at each rational point (cf. \cite{BDBLS05}). See \cite{LZ19} for a similar property for $A_{\chi_4, sq}$. This is one manifestation of Paley-Wiener's theory.

Using the Mellin transforms described in the previous section and the Mellin isometry, we then collect the following relations:
\begin{lemm}[Mellin isometries] \label{mellin-iso}
For all $n,m\geq1$,
\begin{eqnarray*}
\frac{1}{2\pi}\int_{-\infty}^\infty \left|\zeta\left(\tfrac{1}{2}+it\right)\right|^2 \left(\frac{m}{n}\right)^{it} \frac{dt}{t^2+\tfrac{1}{4}} 
        & = & \sqrt{mn} \int_0^\infty \left\{\frac{1}{nx}\right\}\left\{\frac{1}{mx}\right\} dx, \\
\frac{1}{2}\int_{-\infty}^\infty \left|\zeta\left(\tfrac{1}{2}+it\right)\right|^2 \left(\frac{m}{n}\right)^{it} \frac{dt}{\cosh(\pi t)} 
        & = & \sqrt{mn} \int_0^\infty \left(\frac{1}{mx}-\frac{1}{e^{mx}-1}\right)\left(\frac{1}{nx}-\frac{1}{e^{nx}-1}\right) dx, \\
\frac{1}{2\pi}\int_{-\infty}^\infty \left|L\left(\tfrac{1}{2}+it,\chi_4\right)\right|^2 \left(\frac{m}{n}\right)^{it} \frac{dt}{t^2+\tfrac{1}{4}} 
        & = & \sqrt{mn} \int_0^\infty S\left(\frac{1}{nx}\right) S\left(\frac{1}{m x}\right) dx, \\
2\int_{-\infty}^\infty \left|L\left(\tfrac{1}{2}+it,\chi_4\right)\right|^2 \left(\frac{m}{n}\right)^{it} \frac{dt}{\cosh(\pi t)} 
        & = & \sqrt{mn} \int_0^\infty \frac{dx}{\cosh(nx)\cosh(m x)}.     
\end{eqnarray*}
\end{lemm}
The closed-form identities in the next section are all obtained from the right hand side of the above expressions, as in \cite{BC13a} for the first one. 

The inner product $\int_0^\infty \left\{\frac{1}{nx}\right\}\left\{\frac{1}{mx}\right\} dx$ shows up in B\'aez-Duarte's criterion for RH. The other r.h.s. can be interpreted as inner products in some Nyman-Beurling-B\'aez-Duarte type criteria (cf. \cite{DH21b,LZ19}). See also \cite{dR07,DFMR13} for other generalizations.

It is a remarkable feature of $L(\cdot,\chi_4)$ that the function $\cosh$ appears both on the left and right hand sides of the last identity. See also \cite[Sec. 5.2]{BPY01} for interesting approximations of $L(\cdot,\chi_4)$ through probability representations.

\section{Closed-form expressions for the twisted mean squares of $\zeta$ and $L(\cdot,\chi_4)$} \label{review}


From the previous expression of $I_{\zeta, sq}(m/n)$ and the original 
Vasyunin's formula, one obtains (see \cite{BC13a} p. 225):

\begin{theo}[\cite{Vas95}] \label{th:Vasyunin} For coprime $m,n\geq 1$, with the convention $\sum_{\varnothing}=0$, we have
\begin{multline*}
\sqrt{nm} \int_{-\infty}^\infty \left|\zeta\left(\tfrac{1}{2}+it\right)\right|^2 \left(\frac{m}{n}\right)^{it} \frac{dt}{t^2+\tfrac{1}{4}}     
    =\  \pi(\log 2\pi -\gamma)(n+m) + \pi(m-n)\log\left(\frac{n}{m}\right) \\
       - \pi^2 \sum_{k=1}^{n-1}\left\{\frac{k m}{n}\right\} \cot\left(\frac{k \pi}{n}\right) - \pi^2 \sum_{l=1}^{m-1} \left\{\frac{l n}{m}\right\} \cot\left(\frac{l \pi}{m}\right).
\end{multline*}
\end{theo}
See also, e.g., \cite[Sec. 8]{BDBLS00} for non coprime $n,m$.

Bettin and Conrey showed in~\cite{BC13b} a remarkable property of cotangent sums: the function
\bean 
c(x)=-\sum_{a=1}^{k-1} \frac{a}{k} \cot \left(\frac{\pi a h}{k} \right),\ \ x=h/k, \ k>0,\  {\rm gcd}(h,k)=1,
\eean 
satisfies a so-called {\it reciprocity formula}: 
\begin{equation*}
    xc(x)+c\left(\frac{1}{x}\right)-\frac{1}{\pi k} = g(x),
\end{equation*} 
in the vein of Dedekind's reciprocity law. Strikingly, 
the function $g$ is smooth and can be written as $\d g(x)=i\frac{x}{2}\psi(x)$ for $x>0$, where $\psi$ is the period function defined by 
\begin{eqnarray*}
\psi(z) = E(z) - \frac{1}{z}\ E\left(-\frac{1}{z}\right),\quad
E(z) = 1- 4 \sum_{n=1}^\infty d(n) e^{2\pi i nz},
\end{eqnarray*}
for $\Im z>0$ initially, and by analytic continuation to $\C\setminus \R_{\leq0}$.
This remarkable property is also connected to the fundamental relation $A_{\zeta, ch}(z)=\frac{i\pi}{4}\psi(z)+r(z)$, where $r$ is a smooth term due to the polar structure of $\zeta$, see \cite{DH24,DN24}.

The following formula can then be obtained from the period function $\psi$, or from the autocorrelation function $A_{\zeta,ch}$ by elementary means as in \cite{DH21a}:

\begin{theo}[\cite{BC13b,DH21a}] \label{DH21}
For coprime $m,n\geq 1$, with the convention $\sum_{\varnothing}=0$, we have
\begin{multline*}
\sqrt{nm} \int_{-\infty}^\infty \left|\zeta\left(\tfrac{1}{2}+it\right)\right|^2 \left(\frac{m}{n}\right)^{it} \frac{dt}{\cosh(\pi t)}     =\    - 1 + (\log 2\pi -\gamma) (n+m) + (m-n)\log\left(\frac{n}{m}\right) \\
     - \pi \sum_{k=1}^{n-1}\frac{m k}{n} \cot\left(\frac{m k}{n}\pi \right) - \pi \sum_{l=1}^{m-1} \frac{n l}{m} \cot\left(\frac{n l}{m} \pi \right).
\end{multline*}
\end{theo}
See also \cite{DH21a} for non coprime $n,m$ with an abuse of notation. We refer to \cite{Bet15,MR16,ABB17} and references therein for many aspects and generalizations of cotangent sums.

In \cite{LZ19} Lewis and Zagier study a criterion for GRH, especially for $L(\cdot,\chi_4)$. This involves a Gram matrix with a coefficient $c_{m,n}$, which turns out to be a closed form of cotangent sum, related to quantum modular forms, see \cite{Zag10, LZ19, BFOR17}. 

In \cite[Prop 1-2, p.6]{LZ19}, the authors provide a formula for 
$$ c_{n,m}:=\frac{4}{\pi} \int_0^\infty S(nt)S(mt)\frac{dt}{t^2} = \frac{4}{\pi}  mn \int_0^\infty S\left(\frac{1}{nx}\right) S\left(\frac{1}{m x}\right) dx.$$
Their proof is particularly short and elegant, based on Euler's formula for $\pi^2/\sin^2(\pi x)$. 
With Mellin isometry, their formula reads as the following 

\begin{theo}[\cite{LZ19}] \label{LZ19}
Let $m,n\in\N$, and set for $j,k\in\N$: 
$$ a_{m,n}^{j,k}=\max\left(\frac{4j+1}{4m},\frac{4k+1}{4n}\right), \quad b_{m,n}^{j,k}=\min\left(\frac{4j+3}{4m},\frac{4k+3}{4n},\frac{1}{2}\right).$$
Therefore:
\begin{eqnarray*}
\sqrt{nm} \int_{-\infty}^\infty \left|L\left(\tfrac{1}{2}+it,\chi_4\right)\right|^2 \left(\frac{m}{n}\right)^{it} \frac{dt}{t^2+\tfrac{1}{4}} 
    \ =\   \frac{\pi^2}{2} \sum_{\substack{0\leq j\leq m/2 \\ 0\leq k\leq n/2}} \left( \cot\left(a_{m,n}^{j,k}\pi\right) -\cot\left(b_{m,n}^{j,k}\pi\right) \right) \mathds{1}_{a_{m,n}^{j,k}<b_{m,n}^{j,k}}. 
\end{eqnarray*}
\end{theo}
\smallskip

It is then tempting to get the corresponding formula for the $ch-$weight, thus completing the "rectangle" in the introduction. This is the content of the following
\begin{theo} \label{new}
Let $\gcd(n,m)=1$. If $n+m$ is odd, then
\begin{eqnarray*}
\sqrt{nm} \int_{-\infty}^\infty \left|L\left(\tfrac{1}{2}+it,\chi_4\right)\right|^2 \left(\frac{m}{n}\right)^{it} \frac{dt}{\cosh(\pi t)} 
    & = & \frac{\pi}{4} m \sum_{k=0}^{n-1} \frac{(-1)^k}{\cos\!\left(\frac{1+2k}{2n}m\pi\right)} +
    \frac{\pi}{4} n \sum_{\ell=0}^{m-1}  \frac{(-1)^\ell}{\cos\!\left(\frac{1+2\ell}{2m}n\pi\right)}. 
\end{eqnarray*}
If $n+m$ is even, i.e. $n$ and $m$ are odd, then
\begin{multline*}
2\sqrt{nm} \int_{-\infty}^\infty \left|L\left(\tfrac{1}{2}+it,\chi_4\right)\right|^2 \left(\frac{m}{n}\right)^{it} \frac{dt}{\cosh(\pi t)} 
    \ =\  - \pi\frac{m}{n} \sum_{\substack{k\leq n-1 \\ k\neq \frac{n-1}{2}}} \frac{(-1)^k\ k}{\cos\!\left(\frac{1+2k}{2n}m\pi\right)} \\
     \ \ - \pi\frac{n}{m} \sum_{\substack{\ell\leq m-1 \\ \ell\neq \frac{m-1}{2}}} \frac{(-1)^\ell\ \ell}{\cos\!\left(\frac{1+2\ell}{2m}n\pi\right)} -(-1)^{\frac{n+m}{2}}. 
\end{multline*}
\end{theo}

The function $1/\cos$ appears naturally in our proof, but can be written as:
$$\frac{2}{\cos \theta}=\cot\left(\frac{\pi}{4}+\frac{\theta}{2}\right) + \cot\left(\frac{\pi}{4}-\frac{\theta}{2}\right).$$
Hence the previous integrals are also cotangent sums. For instance, when $n+m$ is odd, we have 
$$\sum_{k=0}^{n-1} \frac{(-1)^k}{\cos\!\left(\frac{1+2k}{2n}m\pi\right)} =  \sum_{k=0}^{n-1} (-1)^k \cot\left(\frac{\pi}{4}+\frac{1+2k}{4n}m\pi \right),$$
by reindexing $k\to n-1-k$ in the second $\sum_{k=0}^{n-1}\cot$.
We now prove Theorem \ref{new}.

\section{Proof of the new formula}
\label{proof}

The closed-form expression in Theorem \ref{new} stems from the following

\begin{prop} \label{closed-new}
Let $\gcd(n,m)=1$. If $n+m$ is odd, then
    \begin{eqnarray*}
    \int_{0}^\infty \frac{dx}{\cosh(nx)\cosh(mx)} & = & \frac{\pi}{2n} \sum_{k=0}^{n-1} \frac{(-1)^k}{\cos\!\left(\frac{1+2k}{2n}m\pi\right)} +
    \frac{\pi}{2m} \sum_{\ell=0}^{m-1}  \frac{(-1)^\ell}{\cos\!\left(\frac{1+2\ell}{2m}n\pi\right)}. 
\end{eqnarray*}

If $n+m$ is even, then
\begin{eqnarray*}
    \int_{0}^\infty \frac{dx}{\cosh(nx)\cosh(mx)} & = & - \frac{\pi}{n^2} \sum_{\substack{k\leq n-1 \\ k\neq \frac{n-1}{2}}} \frac{(-1)^k\ k}{\cos\!\left(\frac{1+2k}{2n}m\pi\right)} - 
    \frac{\pi}{m^2} \sum_{\substack{\ell\leq m-1 \\ \ell\neq \frac{m-1}{2}}} \frac{(-1)^\ell\ \ell}{\cos\!\left(\frac{1+2\ell}{2m}n\pi\right)} -\frac{(-1)^{\frac{n+m}{2}}}{nm}. 
\end{eqnarray*}
\end{prop}


\subsection{First case}

Let us prove Proposition \ref{closed-new} in the case where $m,n \geq 1$ are coprime and $n+m$ is odd. 
We consider the function 
$$f: z \mapsto \frac{1}{\cosh(n z)\cosh(m z)},$$ 
which is meromorphic on $\C$ with poles at $$z_{k,n} = \frac{2k+1}{2n} i \pi \ , \ k \in \ZZ \ \ \textrm{and} \  z_{\ell,m} = \frac{2\ell+1}{2m} i \pi \ , \ \ell \in \ZZ. $$ 
According to the assumptions on $n$ and $m$, we have $z_{k,n} \neq z_{l,m}$ for any $k,\ell \in \ZZ$, which means that the poles of $f$ are all simple. Indeed, they are double if and only if: 
$$ \frac{2k+1}{2n} =\frac{2\ell+1}{2m} \iff
m(2k+1) = n(2\ell+1) \iff
m-n = 2(\ell n -km).
$$

\medskip

Consider a pole at $z_{k,n}$ for some $k \in \ZZ$. The residue at $z_{k,n}$ is given by:
\begin{eqnarray*}
\text{Res}(f, z_{k,n}) &=& \lim_{z \to z_{k,n}} (z - z_{k,n}) \frac{1}{\cosh(mz)\cosh(nz)} \\
&=& \frac{1}{\cosh(mz_{k,n})} \lim_{z \to z_{k,n}} \frac{z - z_{k,n}}{\cosh(nz)} \\
&=& \frac{1}{\cosh(m z_{k,n})\ n \sinh(n z_{k,n})}.
\end{eqnarray*}

Furthermore we have: 
\begin{equation*}
\sinh(nz_{k,n}) = \sinh\left(i \frac{2k+1}{2} \pi\right) = i \sin\left(\frac{2k+1}{2} \pi\right) = i (-1)^k,
\end{equation*}
and then:

\begin{equation*}
\text{Res}(f, z_{k,n}) = -\frac{i}{n} \frac{(-1)^k}{\cos\left(\left(k+\frac{1}{2}\right) \frac{m}{n} \pi \right)}.
\end{equation*}

Similarly, we compute, for each $\ell \in \mathbb{Z}$:
\begin{equation*}
\text{Res}(f, z_{\ell,m}) = -\frac{i}{m} \frac{(-1)^l}{\cos\left(\left(l+\frac{1}{2}\right) \frac{n}{m} \pi \right)}.
\end{equation*}

For $R > 0$, we consider the contour:
\begin{equation*}
\Gamma_R : -R \rightarrow R \rightarrow R+i \pi \rightarrow -R+i\pi \rightarrow -R.
\end{equation*}

Due to the parity of $n$ and $m$, we have $f(z+i\pi) = -f(z)$, from which we deduce that: 
\begin{equation*}
 \int_{0}^\infty \frac{dx}{\cosh(nx)\cosh(mx)} = \frac{1}{4} \lim_{R \to \infty} \int_{\Gamma_R} f(z) dz.
\end{equation*}

The poles of $f$ contained in $\Gamma_R$ are the $z_{k,n}$ for $0 \leq k \leq n-1$ and the $z_{l,m}$ for $0 \leq \ell \leq m-1$. The residue theorem then gives:
\begin{eqnarray*}
4 \int_{0}^\infty \frac{dx}{\cosh(nx)\cosh(mx)} &=& 2 i \pi\left ( \sum_{k=0}^{n-1} \text{Res}(f,z_{k,n}) + \sum_{\ell=0}^{m-1} \text{Res}(f,z_{\ell,m}) \right) \\
&=& \frac{2 \pi}{n} \sum_{k=0}^{n-1} \frac{(-1)^k}{\cos\!\left(\frac{1+2k}{2n}m\pi\right)} +
    \frac{2 \pi}{m} \sum_{\ell=0}^{m-1}  \frac{(-1)^\ell}{\cos\!\left(\frac{1+2\ell}{2m}n\pi\right)}. 
\end{eqnarray*}

\subsection{Second case}

Let us now study the case where $m,n \geq 1$ are coprime and both odd. We need to figure out some tricks. The following might be the simplest, considering the function: $$g : z \mapsto \frac{z}{\cosh(n z)\cosh(m z)}.$$

Due to the parity of $n$ and $m$, we deduce that: 
\begin{eqnarray*}
 \lim_{R \to \infty} \int_{\Gamma_R} g(z) dz &=& \lim_{R \to \infty} \left[\int_{-R}^R g(x) dx -  \int_{-R}^R g(x+i\pi) dx \right] \\
 &=& -\lim_{R \to \infty} \int_{-R}^R \frac{x+i\pi}{\cosh(nx) \cosh(mx)} dx \\
 &=& -i\pi \lim_{R \to \infty} \int_{-R}^R \frac{dx}{\cosh(nx) \cosh(mx)} \\
 &=& -2i\pi \int_{0}^\infty \frac{dx}{\cosh(nx)\cosh(mx)}.
\end{eqnarray*}

The poles of $g$ within $\Gamma_R$ are the same as the poles of $f$, but in this case there exists a unique double pole at $$z_{k,\frac{n-1}{2}} = z_{\ell,\frac{m-1}{2}} = i \frac{\pi}{2}.$$

For the simple poles we have:
\begin{equation*}
\text{Res}(g, z_{k,n})  =  \frac{-i(-1)^k}{n \cos\left(\frac{m(1 + 2k)}{2n}\pi\right)} z_{k,n} \ , \quad 
\text{Res}(g, z_{\ell,m})  =  \frac{-i(-1)^\ell}{n \cos\left(\frac{n(1 + 2\ell)}{2m}\pi\right)} z_{\ell,m}.
\end{equation*}

For the double pole we have:
\begin{equation*}
\text{Res}\left(g, i\frac{\pi}{2}\right) = (-1)^{\frac{n+m}{2}+1} i \frac{\pi}{2mn}.
\end{equation*}

So far, the residue theorem gives:
\begin{eqnarray*}
    \int_{0}^\infty \frac{dx}{\cosh(nx)\cosh(mx)} &=& - \sum_{\substack{k\leq n-1 \\ k\neq \frac{n-1}{2}}} \frac{-i(-1)^k}{n \cos\left(\frac{m(1 + 2k)}{2n}\pi\right)} z_{k,n} - 
     \sum_{\substack{\ell\leq m-1 \\ \ell\neq \frac{m-1}{2}}} \frac{-i(-1)^\ell}{n \cos\left(\frac{n(1 + 2\ell)}{2m}\pi\right)} z_{\ell,m} -\frac{(-1)^{\frac{n+m}{2}}}{nm}\\
     &=& -\frac{\pi}{2n^2} \sum_{\substack{k\leq n-1 \\ k\neq \frac{n-1}{2}}} \frac{(-1)^k (2k+1)}{\cos\left(\frac{m(1 + 2k)}{2n}\pi\right)} -\frac{\pi}{2m^2} \sum_{\substack{\ell\leq m-1 \\ \ell \neq \frac{m-1}{2}}} \frac{(-1)^\ell (2\ell+1)}{\cos\left(\frac{n(1 + 2\ell)}{2m}\pi\right)}  -\frac{(-1)^{\frac{n+m}{2}}}{nm}.
\end{eqnarray*}

We now notice, since $m$ is odd, that
\begin{eqnarray*}
\cos\left(m \frac{2(n-k-1)+1}{2n} \pi \right) = - \cos\left(m \frac{2k+1}{2n} \pi \right),
\end{eqnarray*}
from which we deduce, since $n$ is odd: 
$$\frac{(-1)^{n-1-k}}{\cos\left(m \frac{2(n-k-1)+1}{2n} \pi \right) }+\frac{(-1)^k}{\cos\left(m \frac{2k+1}{2n} \pi \right)} = 0,$$
and finally: 
$$\sum_{\substack{k\leq n-1 \\ k\neq \frac{n-1}{2}}}\frac{(-1)^k}{\cos\left(\frac{m(1 + 2k)}{2n}\pi\right)} = 0.$$

Reasoning similarly on the second sum allows to simplify the above equation:
\begin{equation*}
\int_{0}^\infty \frac{dx}{\cosh(nx)\cosh(mx)} = -\frac{\pi}{n^2} \sum_{\substack{k\leq n-1 \\ k\neq \frac{n-1}{2}}} \frac{(-1)^k k}{\cos\left(\frac{m(1 + 2k)}{2n}\pi\right)} -\frac{\pi}{m^2} \sum_{\substack{\ell\leq m-1 \\ \ell \neq \frac{m-1}{2}}} \frac{(-1)^\ell \ell}{\cos\left(\frac{n(1 + 2\ell)}{2m}\pi\right)}  -\frac{(-1)^{\frac{n+m}{2}}}{nm},
\end{equation*}
as claimed.
\qed

\section{Moments of a weighted mean square of $\zeta$ and Dirichlet $L-$functions} \label{review-moment}

In this section, we recall and write in a consistent manner the closed forms related to 
\begin{eqnarray*}
m_N(\zeta) & = & \int_{-\infty}^\infty t^{N} \left|\zeta\left(\tfrac{1}{2}+it\right)\right|^2 \frac{dt}{\cosh(\pi t)}\\ 
m_N(\chi) & = & \int_{-\infty}^\infty t^N \left|L\left(\tfrac{1}{2}+it, \chi\right)\right|^{2} \frac{dt}{\cosh(\pi t)}, 
\end{eqnarray*}
obtained resp. in \cite{DH24} and \cite{DRR26}.
These moments are related to a measure on the critical line involving $\zeta$ or $L(\cdot,\chi)$, and do not refer to the so-called moments of the Riemann zeta function or $L-$functions. The latter, which constitute a huge topic, appear as $\int_0^T |\zeta(\tfrac{1}{2}+it)|^{2k} dt$ with $T>0$ ($k$ being an integer or not), or with a weight on $\R$ (see, e.g., \cite[Chap. VII]{Tit86}).

Let us recall the (Hamburger) moment problem in measure theory. Given a real sequence $(m_n)_{n\geq 0}$, does there exist a positive Borel measure $\nu$ such that for all $n\geq0$, $m_n=\int_\R t^n d\nu(t)$? If yes, is $\nu$ unique? In this last case the problem is called determinate. 

The moments $m_N(\zeta)$ and $m_N(\chi)$
constitute determinate problems due to the exponential weight and the polynomial growth of $\zeta$ and $L(\cdot,\chi)$ on the critical line. Rephrased, $t\mapsto \left|\zeta\left(\tfrac{1}{2}+it\right)\right|^2$ is the only non-negative measurable even function $\phi$ on $\R$ for which $\d \int_{-\infty}^\infty t^{2n} \phi(t) \frac{dt}{\cosh(\pi t)} = m_{2n}(\zeta)$ (given below), since $1/\cosh$ does not vanish (on a set in $\R$ of positive measure).

The Gauss sum of a character $\chi$ of conductor $q$ is defined by: $\d\tau(\chi) = \sum _{a=1}^{q} \chi (a) e^{2\pi i\frac{a}{q}}$, which is not zero for a primitive character $\chi$. 
Let $B_j(\cdot)$, resp. $B_j$, denote the $j$th Bernoulli polynomial, resp. the $j$th Bernoulli number, and $\varphi$ the Euler totient function.
The non-central Stirling numbers of the second kind is defined in \cite{Kou82} by the induction: 
$$S_\alpha(n+1,k)=S_\alpha(n,k-1)+(k-\alpha)S_\alpha(n,k),$$ 
with $S_\alpha(n,0)=(-\alpha)^n$ if $n\geq0$, and $S_\alpha(0,k)=0$ if $k\geq1$, where $\alpha\in\R$. The classical Stirling numbers of the second kind $S(n,k)$ corresponds to the case $\alpha=0$.


\begin{theo}[\cite{DRR26}] \label{DRR26}
Let $\chi$ be a primitive character $\mod q\geq3$. 

Then for all integer $N\geq0$, 
\begin{multline*}
\int_{-\infty}^\infty \left|L\left(\tfrac{1}{2}+it, \chi\right)\right|^{2} \frac{t^N\ dt}{\cosh(\pi t)} \   = \  -\chi(-1) (-1)^{\lfloor \frac{N}{2}\rfloor} B_N\left(\tfrac{1}{2}\right) \frac{\varphi(q)}{q} \\
    + \sum_{k=1}^{N+1} S_{-\frac{1}{2}}(N,k-1) \frac{(2 \pi i q)^k}{q^2 k^2} \sum_{1\leq a,b\leq q}  B_{k}\left(\tfrac{a}{q}\right) B_{k}\left(\tfrac{b}{q}\right) \scr T_{\chi,N}(ab)\ e^{2\pi i\tfrac{ab}{q}}\ ,
\end{multline*}
where\  $\ \scr T_{\chi,N}(c) = -i^N \left[\tau(\b \chi)\ \chi(c) +(-1)^N \tau(\chi)\ \b \chi(c) \right]$, $c\in \Z$.
\end{theo}
 
 See \cite{DRR26} for the motivations which include the exact characterization of $|L(\tfrac{1}{2}+it)|$ and its relation with the twisted period function
defined for $\Im(z)>0$ by
\begin{eqnarray*}
\psi_\chi(z) & = & E_\chi(z) - \frac{\chi(-1)}{z}\ E_{\b\chi}\left(-\frac{1}{z}\right),\\
E_\chi(z) & = & \frac{1}{\tau(\chi)} \sum_{n=1}^\infty\chi(n)d(n) e^{2\pi i\tfrac{nz}{q}}. 
\end{eqnarray*}

Let us recall some well known expressions that will be used just below (see, e.g., \cite{Ten22,Ap08,Coh07}):
\begin{eqnarray*}
B_{2n}\left(\frac{1}{2}\right) = \left(\frac{1}{2^{2n-1}} -1\right) B_{2n},\quad 
\zeta(2j) = \frac{(-1)^{j+1}}{2(2j)!} (2\pi)^{2j} B_{2j}.
\end{eqnarray*}

To obtain the following formula for $m_{2N}(\zeta)$, we can use the same method as in \cite{DRR26} to make appear directly $S_{-\frac{1}{2}}(2n,2j-1)$, or start from \cite{DH24} and use a relation between the Stirling numbers $S_{-\frac{1}{2}}$ and $S_0$. But it is more diverting to take $q=1$ in Theorem \ref{DRR26} (not licit) and add the contribution of the pole of $\zeta$ from \cite{DH24}:

\begin{theo}[\cite{DH24}] \label{DH24}
For all integer $n\geq0$, we have
\begin{multline*}
\int_{-\infty}^\infty \left|\zeta\left(\tfrac{1}{2}+it\right)\right|^2 \frac{t^{2n} \ dt}{\cosh(\pi t)}\ =\ (-1)^n \left[ \frac{\log(2\pi)-\gamma -4n}{2^{2n-1}} + \left(1 - \frac{1}{2^{2n-1}} \right) B_{2n}  \right. \\
      \qquad \left. +\  \frac{1}{2} \sum_{j=1}^{n} S_{-\frac{1}{2}}(2n,2j-1) \frac{(-1)^{j+1}}{j^2} B_{2j}^2 \cdot (2 \pi)^{2j} \right].
\end{multline*}
\end{theo}
These quantities appear in a moment Gram matrix related to some generalizations of the Nyman-Beurling criterion, see \cite{ADH22}.

Note that $m_{2n}(\zeta)=P_{2n}(\pi)$ where $P_{2n}\in \R+\Q_{n}[X^2]$, and
\begin{eqnarray*}
    m_0(\zeta) & = & 2 (\log(2\pi)-\gamma) - 1 = 1.521322\ldots
\end{eqnarray*}

Let us finally collect the first numerical values of $m_N(\zeta)$, compared to those of $m_N(\chi_4)$ from \cite{DRR26}. Let us notice that the values of $m_N(\zeta)$ here are not the same as in \cite{DH24} (up to a factor $1/\pi$) due to the choice of our common normalization: $dt/\cosh(\pi t)$ in both cases here.\\

\renewcommand{\arraystretch}{1.3}
\begin{center}
\begin{tabular}{|c|c|c|c|}
    \hline
    $N$ & $m_N(\zeta)$ & $m_N(\chi_4)$  \\
    \hline
    0 & $1.521322\ldots$ & $0.5$         \\
    1 & $0$             & $0$             \\
    2 & $0.189713\ldots$ & $0.178744\ldots$  \\
    3 & $0$             & $0$             \\
    4 & $0.138255\ldots$ & $0.334811\ldots$  \\
    5 & $0$             & $0$              \\
    6 & $0.323446\ldots$ & $1.413929\ldots$  \\
    7 & $0$             & $0$              \\
    8 & $1.784230\ldots$ & $9.787352\ldots$ \\
    9 & $0$             & $0$              \\
    10 & $18.344735\ldots$   & $95.077818\ldots$   \\
    \hline
\end{tabular}
\end{center}

We refer to \cite{DRR26} for several pictures and comments on various phenomena as the character varies.

\section{Perspectives}

The similarity between the various reciprocity formulas, with different weights, suggests the existence of a larger theory in which they can all be interpreted, and which would classify the weights and their regularization effect.

A greater challenge at stake is to further develop the theory of reconstruction of a density from its moments (a very unstable problem numerically) or from a dense set of values of its Fourier transform, in order to retrieve either local information about the density or its asymptotic behavior.


\section{Use of LLMs and Acknowledgments}

We used LLMs (OpenAI, Anthropic, Google, Mistral) to complete our search for various closed forms in the literature in our context, to double-check the formulas as they are stated here, and for a few helpful remarks and corrections.
We are glad that some AI found papers with only a few citations, as the papers of Ivi\'c \cite{Ivi03} and Coffey \cite{Cof11}, which we did not know.

S.D. thanks Joseph Najnudel, Berend Ringeling and Emmanuel Royer for insightful discussions.

\end{document}